\documentclass{article}
\usepackage[utf8]{inputenc}
\usepackage{amsfonts}
\usepackage{amsmath}
\usepackage{amsthm}
\usepackage{fullpage}
\usepackage{graphicx}
\usepackage{caption}
\usepackage{subcaption}
\usepackage{tikz}
\tikzstyle{mybox} = [draw=black, fill=white,  thick,
    rectangle, inner sep=10pt, inner ysep=20pt]
\tikzstyle{mybox} = [draw=black, fill=white,  thick,
    rectangle, inner sep=2pt, inner ysep=2pt]
\usetikzlibrary{arrows}

\newtheorem{theorem}{Theorem}
\title{Newton-Ellipsoid Method and its Polynomiography\thanks{This work was carried out in part during a Summer REU program at DIMACS, where the first author was mentor to the second author.}}

\author{Bahman Kalantari\thanks{Department of Computer Science,
        Rutgers University, New Brunswick, New Jersey. {\tt kalantari@cs.rutgers.edu}}
        ~~and~~ Eric Lee\thanks{Department of Mathematics,
        Berkeley. {\tt eric.lee@berkeley.edu}}.}

\date{}

\begin{document}

\maketitle
\begin{abstract}
 We introduce a new iterative root-finding method for complex polynomials, dubbed {\it Newton-Ellipsoid} method.  It is inspired by the Ellipsoid method, a classical method in optimization, and a property of Newton's Method derived in \cite{kalFTA}, according to which at each complex number a half-space can be found containing a root. Newton-Ellipsoid method combines this property, bounds on zeros, together with the plane-cutting properties of the Ellipsoid Method. We present computational results for several examples, as well as corresponding polynomiography. Polynomiography refers to algorithmic visualization of
 root-finding.  Newton's method is the first member of the infinite family of iterations, the {\it basic family}. We also consider general versions of this ellipsoid approach where Newton's method is replaced by a higher-order member of the family such as Halley's method.
\end{abstract}

{\bf Keywords:}    Polynomials; Newton Method; Halley Method; Ellipsoid Method; Polynomiography

\section{Introduction} Solving a polynomial equation,
\begin{equation}
p(z)=a_nz^n+ \cdots +a_1z+a_0=0,  \quad a_i \in \mathbb{C},
\end{equation}
has played a significant role in the history and development of science and mathematics.  Today the problem is present in every branch of mathematics and science, education and more, see \cite{pan97}, \cite{kalbook}.

Approximation of the roots is the only means for solving a general polynomial equations of degree five or higher. In fact, this happens to be true even for quadratics. Newton's method is the most fundamental iterative method and is usually considered in the context of  computing real roots of real polynomials. Cayley \cite{cay} considered Newton's method for complex polynomials and studied its convergence properties for computing the roots of unity, solutions of $z^n-1=0$.

Given a real or complex polynomial $p(z)$, Newton's method continuously applied to each {\it seed} $z_0 \in \mathbb{C}$ generates a sequence $\{z_i \}_{i=1}^\infty$, called an {\it orbit}.  The {\it basin of attraction} of  a root $\theta$ of  $p(z)$ under Newton's method is the set of all complex numbers whose corresponding orbit converges to this root. It is a classic result that Newton's method locally converges quadratically.  For studies on the mathematical properties of the dynamics of Newton's method and, more generally, rational functions, see \cite{beardon91}, \cite{milnor}, \cite{peit}, and \cite{kalbook}.

Cayley was only able to characterize the basin of attractions for quadratic polynomials. When the roots of a quadratic are distinct, each basin of attraction is the Voronoi region of the root,  the set of all complex numbers that are closer to that root than the other root. The analysis of basin of attractions for cubics and beyond is far more complicated.  Pathological cubic examples for Newton's method can be given where the method fails to converge to a root on a set of positive measure: for instance, Smale's example $z^3-2z+2$.  Deciding if, for a cubic polynomial, a particular orbit converges to a root is also highly complex, see Blum et al. \cite{Blum}, who relates this to the undecidability problem over the reals. McMullen \cite{McMullen87} shows that for almost all polynomials of degree at least $4$,  all iterative Newton-like methods will be non-convergent on a set of positive measure.

In this article we introduce a hybrid method combining Newton's method and the ellipsoid method.  We apply our hybrid method to several examples, give  corresponding polynomiography, and make comparisons with Newton's method.  As mentioned earlier, Polynomiography is an algorithmic visualization in solving polynomial equations. In the next section, first we summarize the basic properties of both Newton and ellipsoid methods.  We then describe the Newton-Ellipsoid method, provide the justification behind the method, and present some polynomiography.  Finally, we consider modifications and generalizations of the Newton-ellipsoid method.

\section{The Newton-Ellipsoid Method}

To describe the Newton-Ellipsoid, first we acquaint the reader with the necessary ingredients.

\subsection{Newton's Method}
Newton’s method is an iterative approximation method for finding roots of a complex polynomial $p(z)$. Given initial starting point $z_0$, called a  {\it seed}, the  \textit{Newton's iterate} and \textit{Newton's direction}
at $z_k$ are defined, respectively  as
$$z_{k+1} = z_k - \frac{p(z_k)}{p'(z_k)}, \quad - \frac{p(z_k)}{p'(z_k)}.$$

\subsection{Ellipsoid Method}
The Ellipsoid method is an iterative method used in convex programming. Its most famous application was in linear programming, by formulating problems as an LP feasibility problem and then solving a system of inequalities.  Linear programs have been proved to be solvable in polynomial time \cite{kha79}. The main ingredient in the method is the following fundamental property:

Given  an ellipsoid $E$ of dimension $k$,  and a half-space $H$ passing through its center, we can generate an ellipsoid $E'$ that contains the half-ellipsoid $H \cap E$, with the further property that  the volume of $E'$ is a fraction of the volume of $E$, specifically
$$Vol(E') \leq \rho_k  \cdot Vol(E), \quad \rho_k = \exp(-1/2(k+1))<1 .$$
Repeating this process, one obtains a sequence of ellipsoids whose volumes shrink to zero.  In this way, the ellipsoid method is applied to find a solution to the system of inequalities within a certain tolerance $\epsilon$, in a finite number of iterations, polynomial in $k$ and $\ln \epsilon^{-1}$.  In particular, for $k=2$, $\rho_2 \approx .85$

Formally, given an $k \times k$  symmetric positive definite matrix $B$, and a point $c \in R^k$, an ellipsoid is defined as
$$E(B,c)= \{ x \in R^k :  (x-c)^T B^{-1} (x-c) \leq 1\}.$$
The center of the ellipsoid is $c$, while its axis and their size is determined by the eigenvector and eigenvalues of the matrix, respectively. Given $a \in R^k$, the half-space through the center of the ellipsoid having normal equal to $a$ is
$$H(a)=\{x: a^T(x -c) \leq 0. \}.$$
The half-ellipsoid is
$$HE(a)= E(B,c) \cap H(a).$$
The ellipsoid of least volume that contains $HE(a)$ is $E'(B',c')$, where
$$B'= \frac{k^2}{k^2-1} \bigg (B- \frac{2(Ba)(Ba)^T}{(k+1)a^TBa} \bigg ), \quad c'=c- \frac{Ba}{(k+1)a^TBa}.$$

\subsection{A Property of the Newton Direction}

We identify a complex number $z=x+iy$ with the point $(x,y)$ in the Euclidean plane.

\begin{theorem} \label{FTA} {\bf (Kalantari \cite{kalFTA})}
Given a point $z_0$ in $\mathbb{C}$, consider the half-space $H(z_0)$, passing through $z_0$ whose normal corresponds to ${p(z_0)}{p'(z_0)}$. Then $p(z)$ has a root in $H(z_0)$. \qed
\end{theorem}

\subsection{Bounds on Zeros of Polynomials}  Computing a priori bounds on zeros of polynomials is a classic problem, see e.g. \cite{mc2005}.  Given a polynomial $p(z)$ we can generate bounds on the modulus of its zeros. Specifically, for each integer  $m \geq 2$  one can generate a bound on the modulus of zeros of $p(z)$, see \cite{kalbound}. Using these bounds we can generate a tight initial rectangular region needed to start the Newton-Ellipsoid algorithm. We describe the first few bounds.

Let  $r_m \in [1/2, 1)$ be the unique positive root of the
polynomial $t^{m-1} +t-1$.   Assume $\theta$ is any root of $p(z)$. For $m=2$, $r_2=0.5$ and we have
$$|\theta| \leq {1 \over {r_2}} \max \biggr \{ \bigg| {1 \over {{a_n}}}{a_{n-k+1}} \bigg| ^{1/(k-1)}:
k=2, \dots, n+1 \biggr \}.$$

For $m=3$, $r_3=0.618034$ and we have
$$|\theta| \leq {1 \over {r_3}} \max \biggr \{
\bigg|  {1 \over {{a_n}^2}}  {det \left(
\begin{array}{ccc}
a_{n-1}&a_{n-k+1}\\
                      a_n&a_{n-k+2} \\
\end{array}
\right )} \bigg |^{1/(k-1)}: k=3, \dots, n+2 \biggr \},  \quad a_{-1}=0.$$
For $m=4$, $r_4=0.682328$ and we have
$$|\theta| \leq {1 \over {r_4}} \max \biggr \{
\bigg| {1 \over {{a_n}^3}} { det \left ( \begin{array}{ccc}
 a_{n-1}&
a_{n-2}&a_{n-k+1}\\
                      a_n&a_{n-1}&a_{n-k+2}\\
                      0&a_n&a_{n-k+3}\\
\end{array}
\right )}
                      \bigg |^{1/(k-1)}: k=4, \dots, n+3 \biggr \}, \quad a_{-1}=a_{-2}=0.$$

\subsection{The Algorithm}
\begin{figure}
\begin{center}
\includegraphics[width=2.6in]{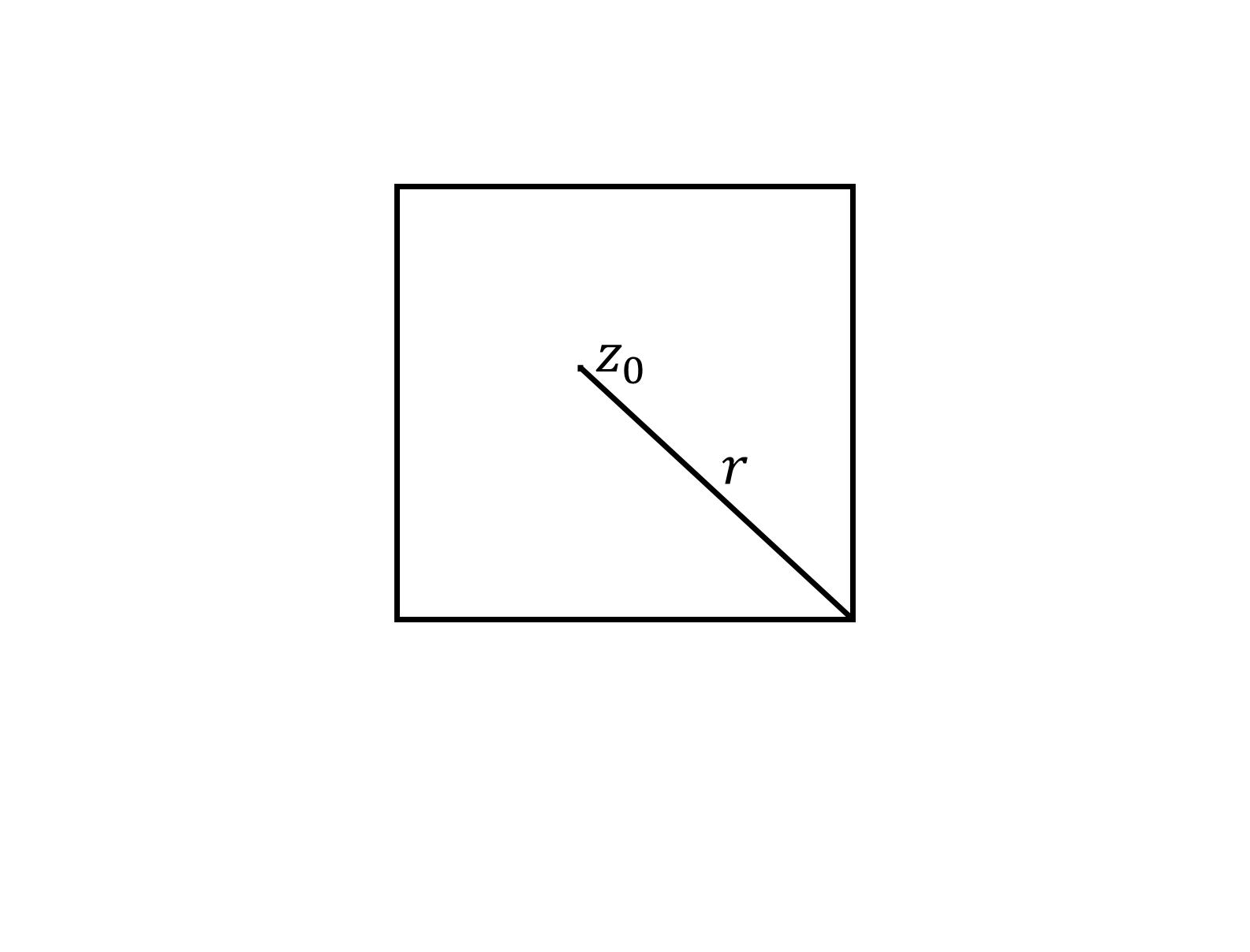}
\includegraphics[width=2.6in]{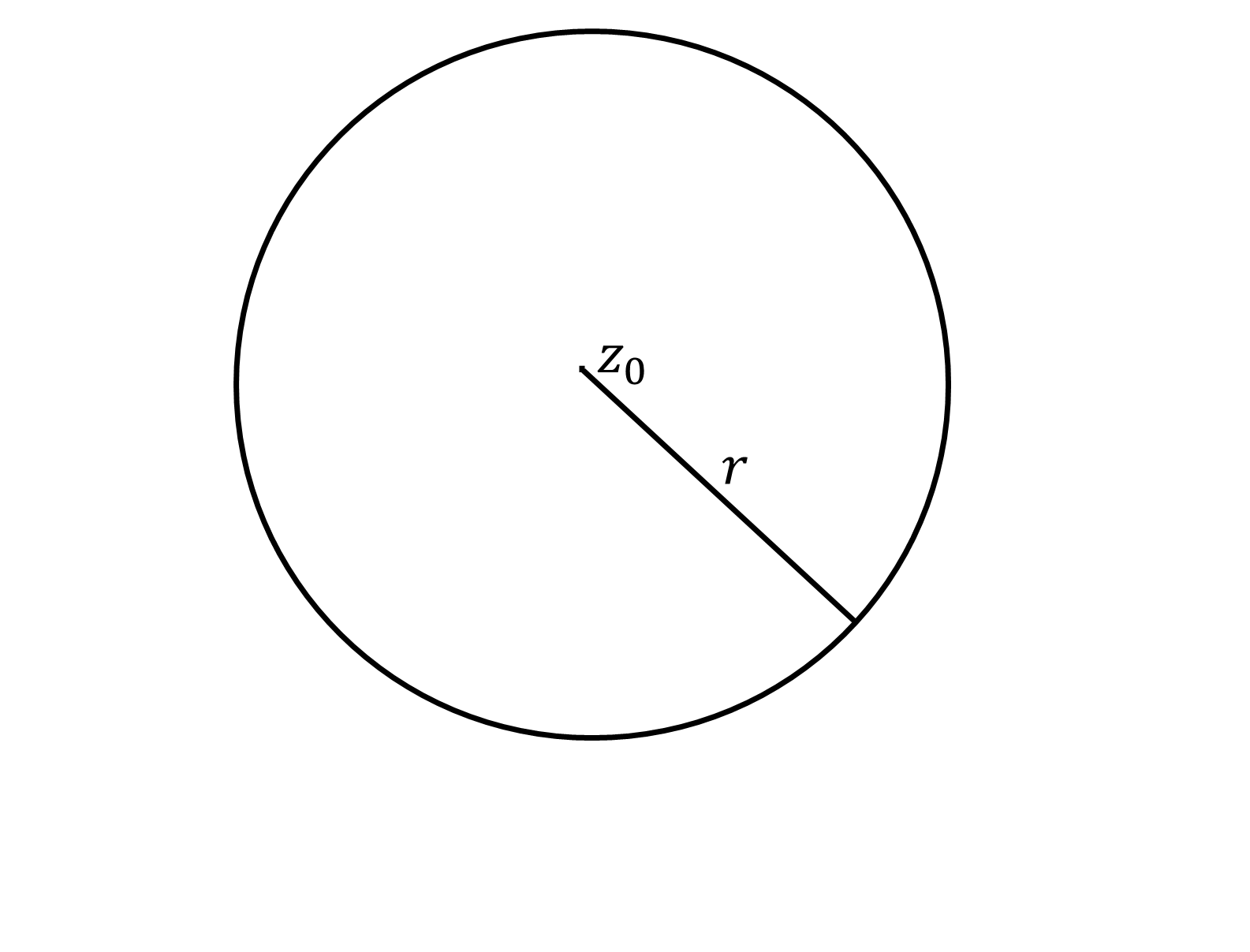}\\
\end{center}
\begin{center}
\caption{An initial box containing the roots (left). A disc containing the roots with center at a given seed $z_0$.} \label{Fig1}
\end{center}
\end{figure}

\begin{figure}
\begin{center}
\includegraphics[width=2.6in]{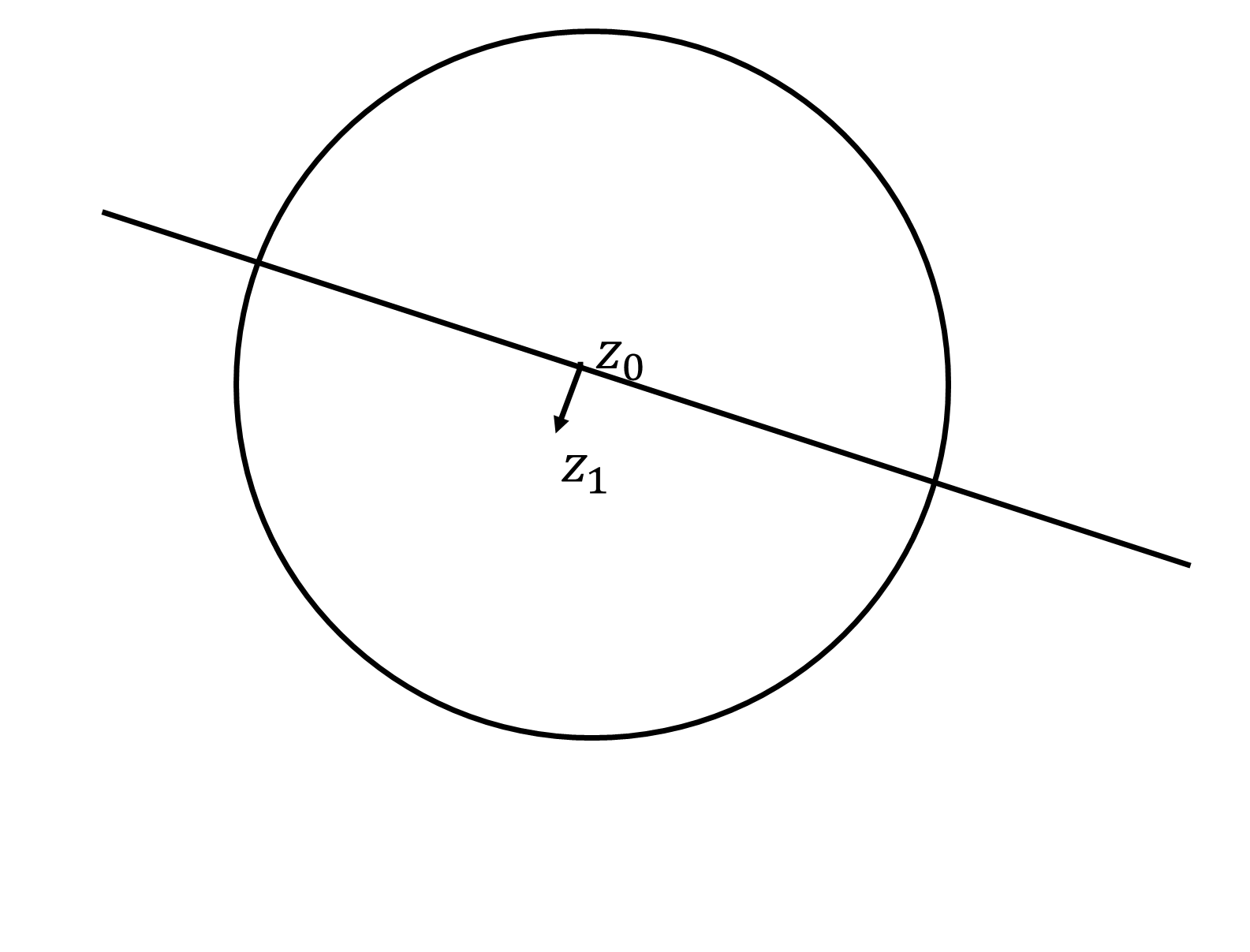}
\includegraphics[width=2.6in]{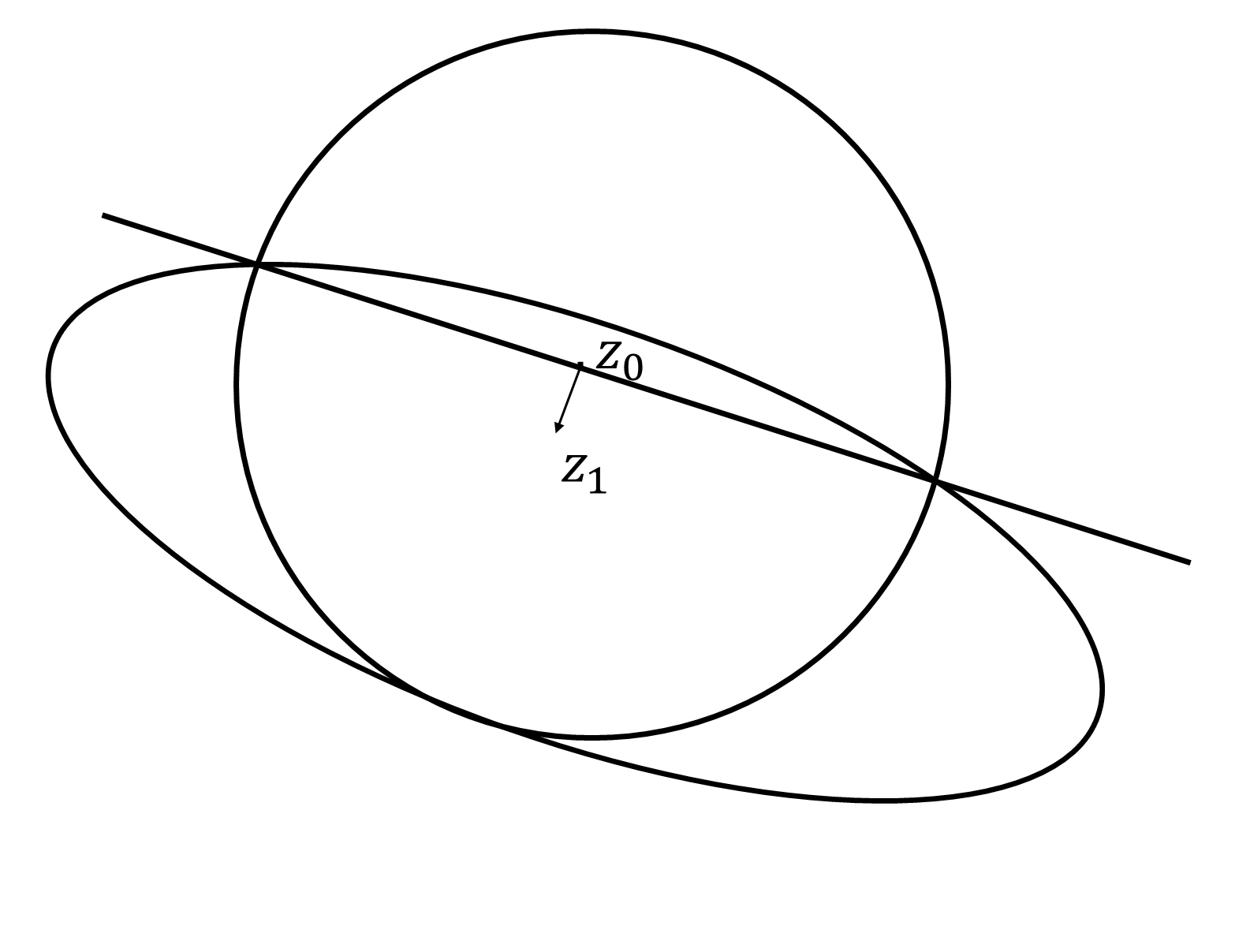}\\
\end{center}
\begin{center}
\caption{Newton's iterate $z_1$ (left) and a smaller ellipsoid determined by Newton's direction.} \label{Fig2}
\end{center}
\end{figure}

\begin{figure}
\begin{center}
\includegraphics[width=2.6in]{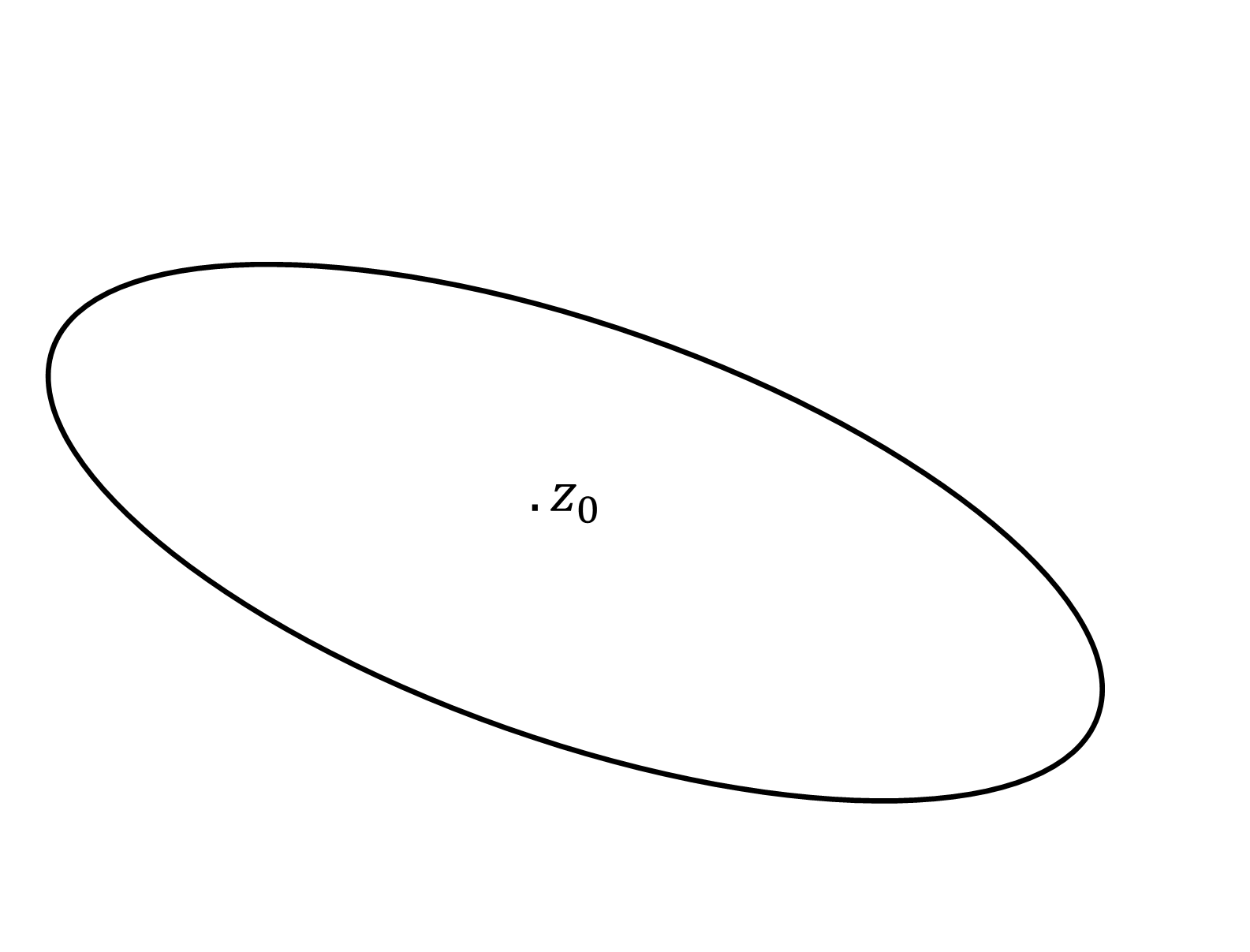}
\includegraphics[width=2.6in]{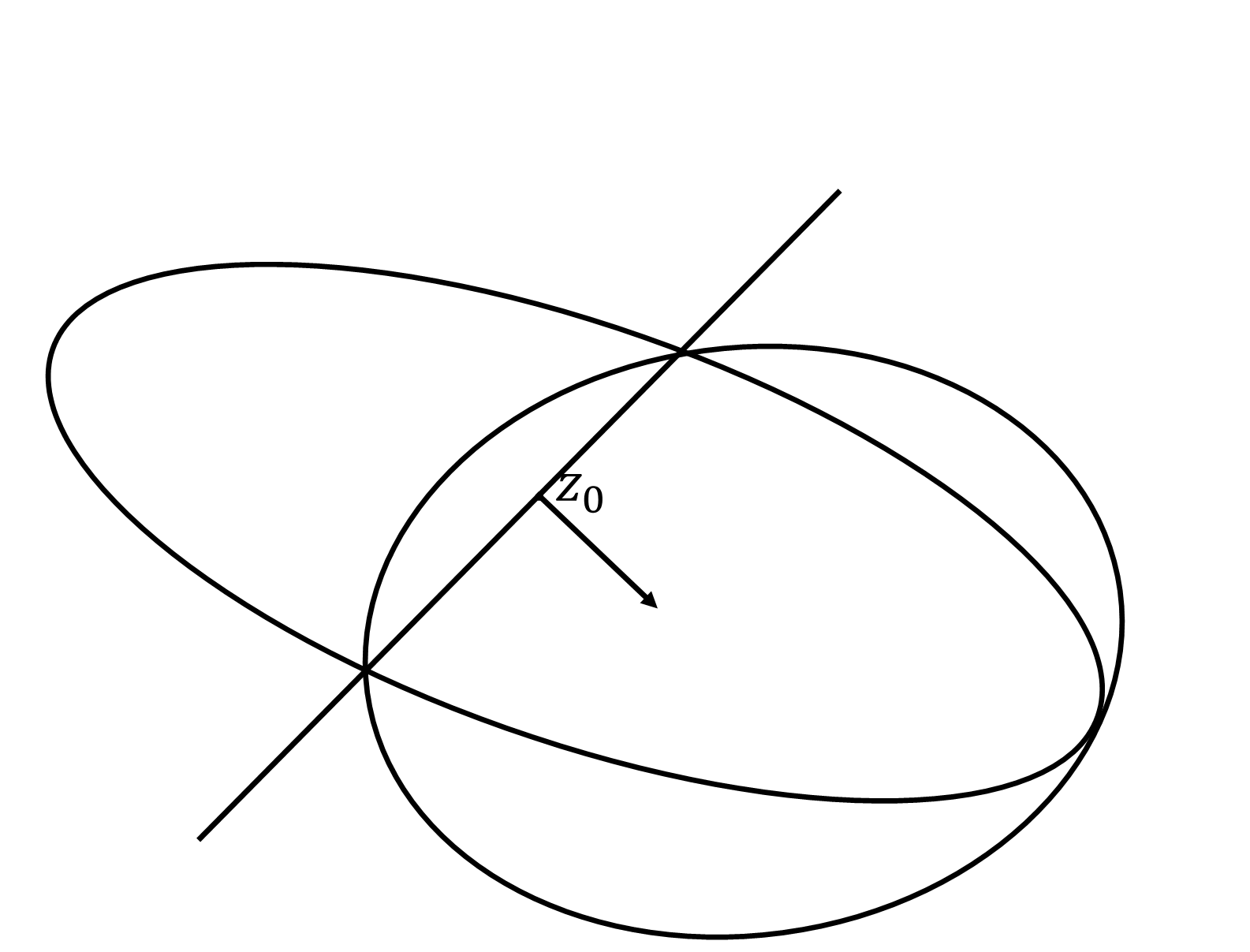}\\
\end{center}
\begin{center}
\caption{The new seed becomes the center of the new ellipsoid and the process repeated.} \label{Fig3}
\end{center}
\end{figure}

We now describe the Newton-Ellipsoid method in further detail. The method is a combination of Newton's method and the Ellipsoid method. The general approach is follows: Given a polynomial $p(z)$, first compute an initial box that contains all the roots. Given an initial seed $z_0$ in the box, we can easily compute the smallest disc that contains the box and has $z_0$ as its center. The disk's radius is simply the distance from $z_0$ to the farthest point in the box (see Figure \ref{Fig1}). This disc is our initial ellipsoid $E_0$.  Next we compute the Newton direction at $z_0$, $-p(z_0)/p'(z_0)$. While Newton's method would take the new point, $z_1$ to be the sum of the old point and its Newton's direction, we instead consider the half space $H(z_0)$ whose normal is $p(z_0)/p'(z_0)$. We draw smaller ellipsoid $E_1$ containing the intersection between the half space and $E_0$. Because by Theorem \ref{FTA} there exists a root in the direction of $-p(z_0)/p'(z_0)$, $E_1$ must contain a root as well. We take the next iterate $z_1$ to be the center of $E_1$ and repeat the process. See Figures \ref{Fig1}-\ref{Fig3}.

The idea is that the Newton-Ellipsoid method - in some instances - will be better-behaved than Newton's Method, and take a larger iterative step. However, while $E_1$ must necessarily contain a root of $p(z)$, $E_2$ may not contain a root, because the Newton direction at the center of $E_1$ may be in a direction of a root that is not in $E_1$. In other words, we run the risk of iterating to a root far away from our original seed. However, with the assumption that Newton direction always points to the nearest root and that the ellipsoid centers will never go outside of the initial box, we can guarantee that the sequence of centers converge to a root of $p(z)$. Formally, the algorithm is described in the box bellow.

\begin{center}
\begin{tikzpicture}
\node [mybox] (box){%
    \begin{minipage}{0.9\textwidth}
{\bf  Newton-Ellipsoid Algorithm} ($p(z)$,  $\epsilon >0$)\

\begin{itemize}

\item  {\bf Step 0.} {\bf (Initialization)} Compute a bound on the modulus of roots. Let $R$ be a rectangle containing the roots. Pick $z_0 \in R$. Let $r$ be the distance from $z_0$ to farthest point of  $R$. Let $E_0(B,z_0)$ be the circle of radius $r$ centered at $z_0$, thus $B=rI$, where $I$ is the $2 \times 2$ identity matrix.

\item {\bf Step 1.}  If $|p(z_0)| < \epsilon$, stop.

\item {\bf Step 2.}  Let $a = \frac{p(z_0)}{p'(z_0)}$, $H(a)= \{z: a^T(z-z_0) \leq 0 \}$.
\item {\bf Step 3.}  Compute the new ellipsoid $E'(B',z')$:

$$z'= z_0 - \frac{Ba}{3\sqrt{a^tBa}},  \quad B'=\frac{4}{3} \bigg (B-\frac{2(Ba)(Ba)^T}{3(a^TBa)} \bigg ).$$

\item {\bf Step 4.}  Set $z_0 = z'$, $B=B'$. Go to Step 1.

\end{itemize}

    \end{minipage}};
\end{tikzpicture}
\end{center}

\subsection{A Justification}
Here we offer a justification that the iterates will converge to a root.
Assume $p(z)$ is a monic complex polynomial of degree $n$. By the fundamental theorem of algebra we can factor $p(z)$ as: $p(z)=(z-\theta_1) \cdots (z-\theta_n)$.
Since $p'(z)= \sum_{i=1}^n p(z)/(z- \theta_i)$,
$$- \frac{p(z)}{p'(z)}= - \frac{1}{\sum_{i=1}^n \frac{1}{(z- \theta_i)}}.$$
We see that when $\theta_j$ is a simple root of $p(z)$ and $z$ is close enough to $\theta_j$ then
$$- \frac{p(z)}{p'(z)} \approx (\theta_j - z).$$
Thus Newton's direction points toward the closest root, namely $\theta_j$.

\subsection{Heuristic Ideas}  If our goal is to compute one or more roots of $p(z)$,  we can pick a random $z_0$ and apply the Newton-Ellipsoid method. If the corresponding sequence does not converge we can pick another initial point and repeat the process. Once a root is approximated by deflation we can repeat the process to find other roots.

\section{Halley-Ellipsoid and Higher Members of the Basic Family}

Given a complex polynomial $p(z)$, the {\it basic family} of iteration functions is the collection of iteration functions defined as
$$B_m(z)=z-p(z) \frac{D_{m-2}(z)} {D_{m-1}(z)}, \quad m=2,3, \dots$$
where  $D_0(z)=1$,  $D_k(z)=0$ for $k <0$, and, $D_m(z)$ satisfies the recurrence relation
$$
D_m(z)= \sum_{i=1}^n (-1)^{i-1}p(z)^{i-1}\frac{p^{(i)}(z)}{i!}D_{m-i}(z).$$

The first two members of the basic family are $B_2(z)$ (Newton's Method) and $B_3(z)$ (Halley's Method).  For the rich  history of history and fundamental properties of the basic family members and applications, see \cite{kalbook}.  For each fixed $m \geq 2$, there exists a disk centered at a root $\theta$ such that for any $z_0$ in this disk, the sequence of fixed point iteration
$z_{k+1}=B_m(z_k)$, $k=0,1,\dots$, is well-defined, and converges to $\theta$.  When  $\theta$ is a simple
root,  the order of convergence is $m$.

There are Voronoi properties in using the basic family members collectively.  For a pointwise convergence property see \cite{kalbook}, and for a uniform convergence proof of, see \cite{KalDCG}. The pointwise convergence property is the following:

\begin{theorem} \label{thm2} For any root $\theta$ of $p(z)$, let $V(\theta)$ be its Voronoi region, i.e. the set of all complex numbers that are closer to $\theta$ than any other root. Then given any $w \in V(\theta)$,  $\lim_{m \to \infty} B_m(w)=\theta.$
\end{theorem}

With respect to the use of the ellipsoid method, we can define the \textit{basic family direction} of order $m$ at a given $z_0$ to be
$$-p(z_0) \frac{D_{m-2}(z_0)} {D_{m-1}(z_0)}.$$

Analogous to Newton-Ellipsoid method, we define the generalized ellipsoid method using the above direction to define a hyperplane passing through each $z_0$.  However, a generalization of Theorem \ref{FTA} needs to be proved. Nonetheless, we will offer present polynomiography for higher order members as well.

\subsection{Polynomiography}
\begin{figure}[h!]
\centering
\includegraphics[width=2.6in]{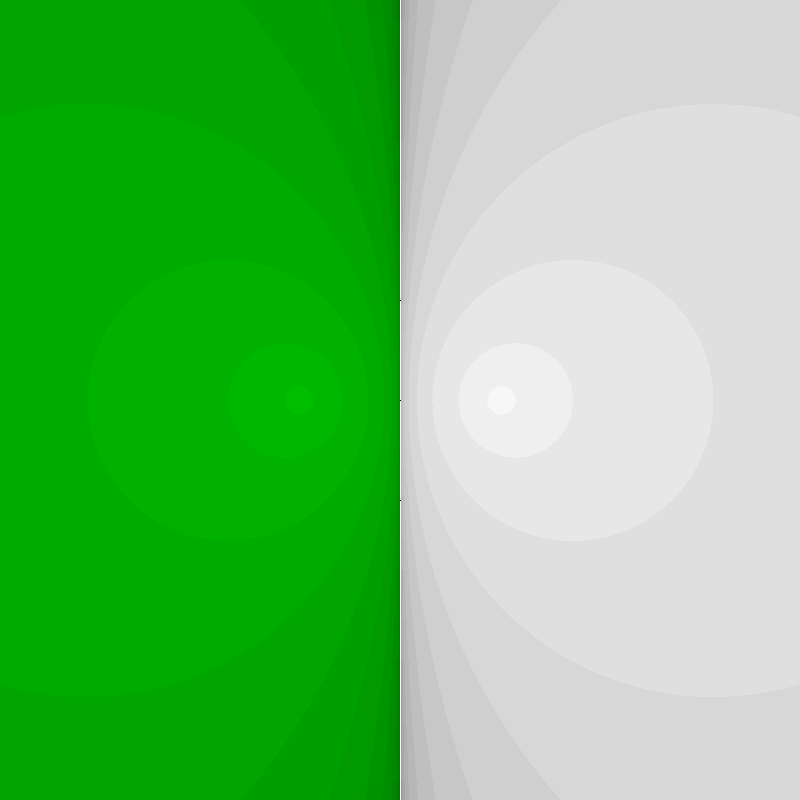}
\includegraphics[width=2.6in]{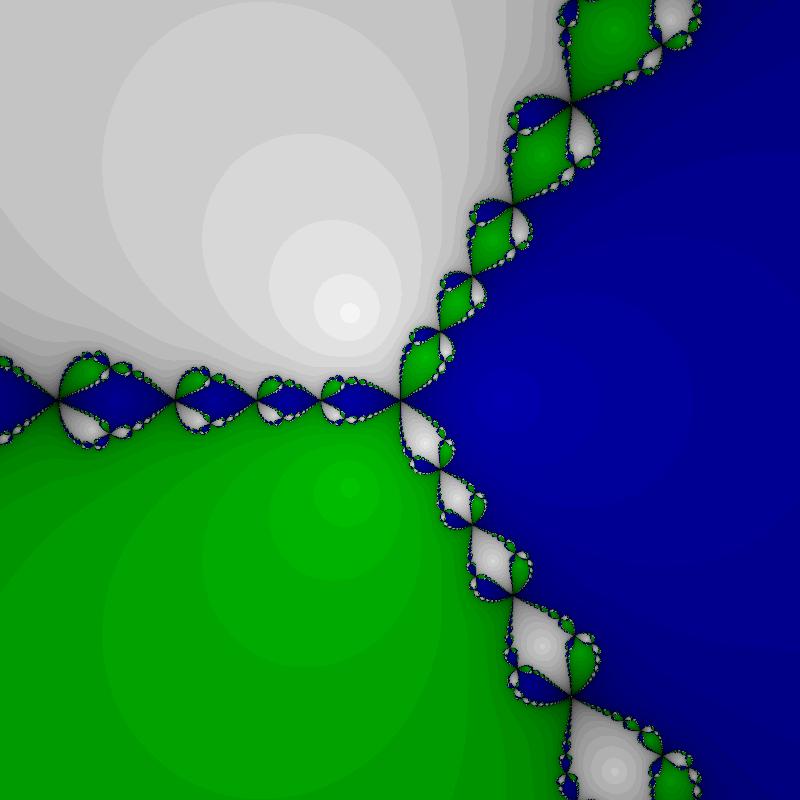}
\caption{Polynomiography of $z^2-1$ (left) and $z^3-1$ under Newton's method.} \label{FigNewton}
\end{figure}

\begin{figure}[h!]
\centering
\includegraphics[width=2.6in]{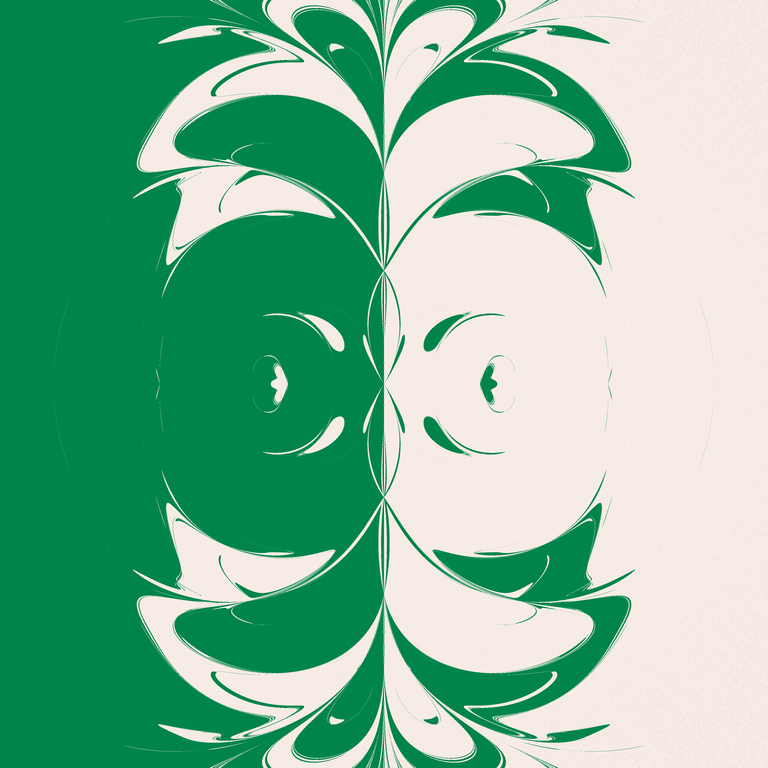}
\includegraphics[width=2.6in]{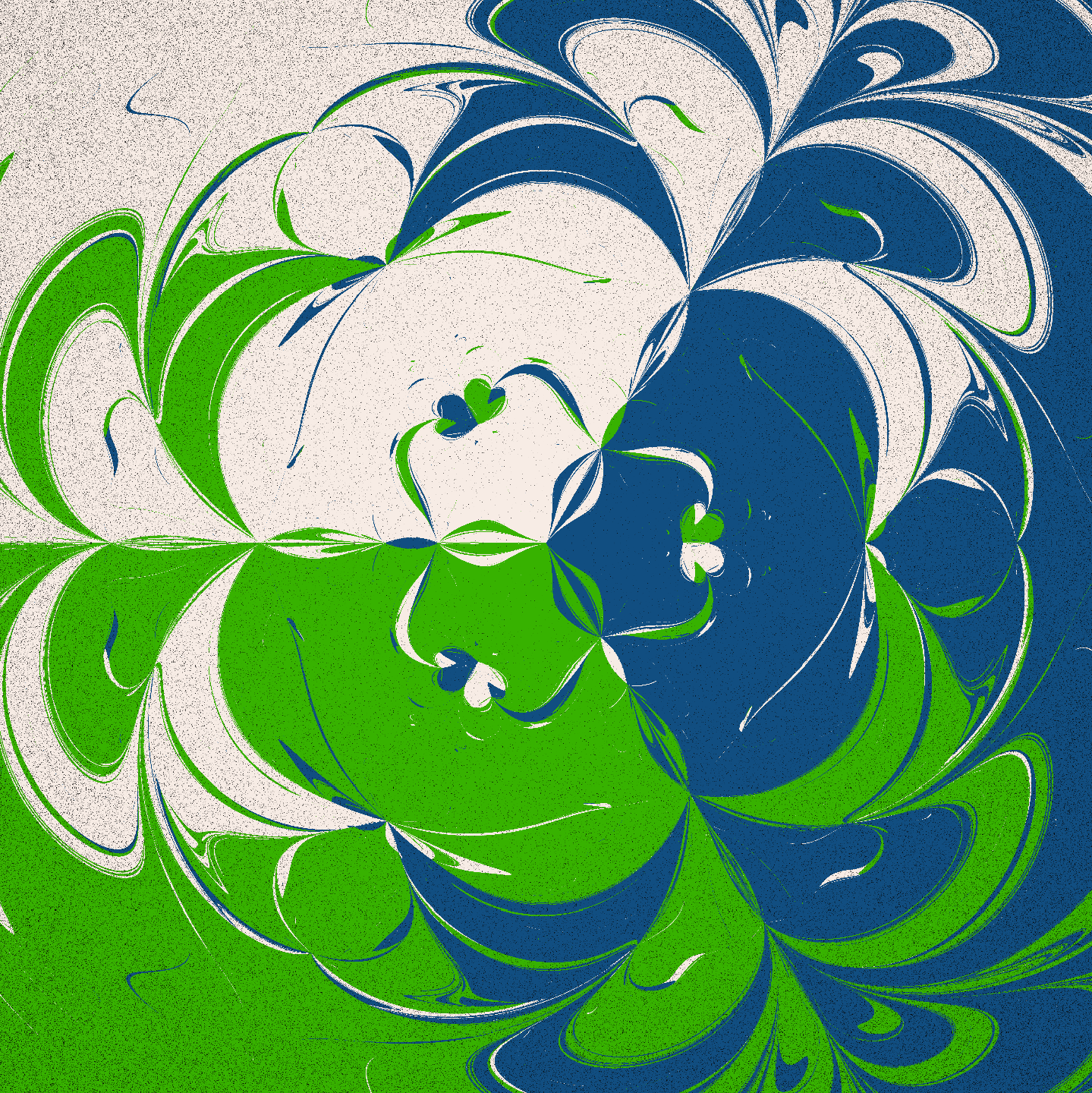}
\caption{Polynomiography of $z^2-1$ (left) and $z^3-1$ under Newton-Ellipsoid method.} \label{Newton-Ellipsoid}
\end{figure}

\begin{figure}[h!]
\centering
\includegraphics[width=2.6in]{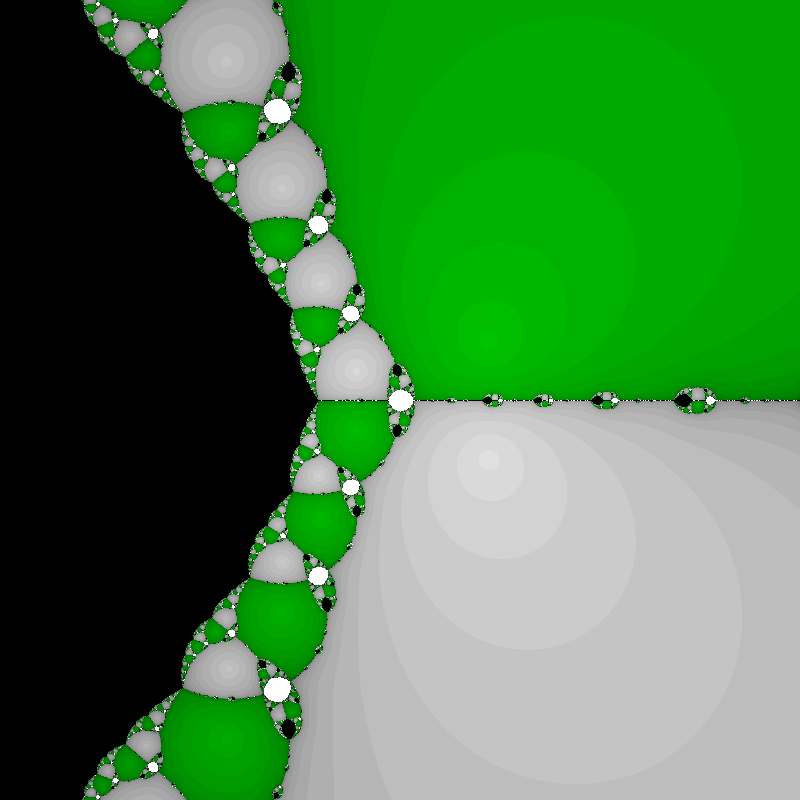}
\includegraphics[width=2.6in]{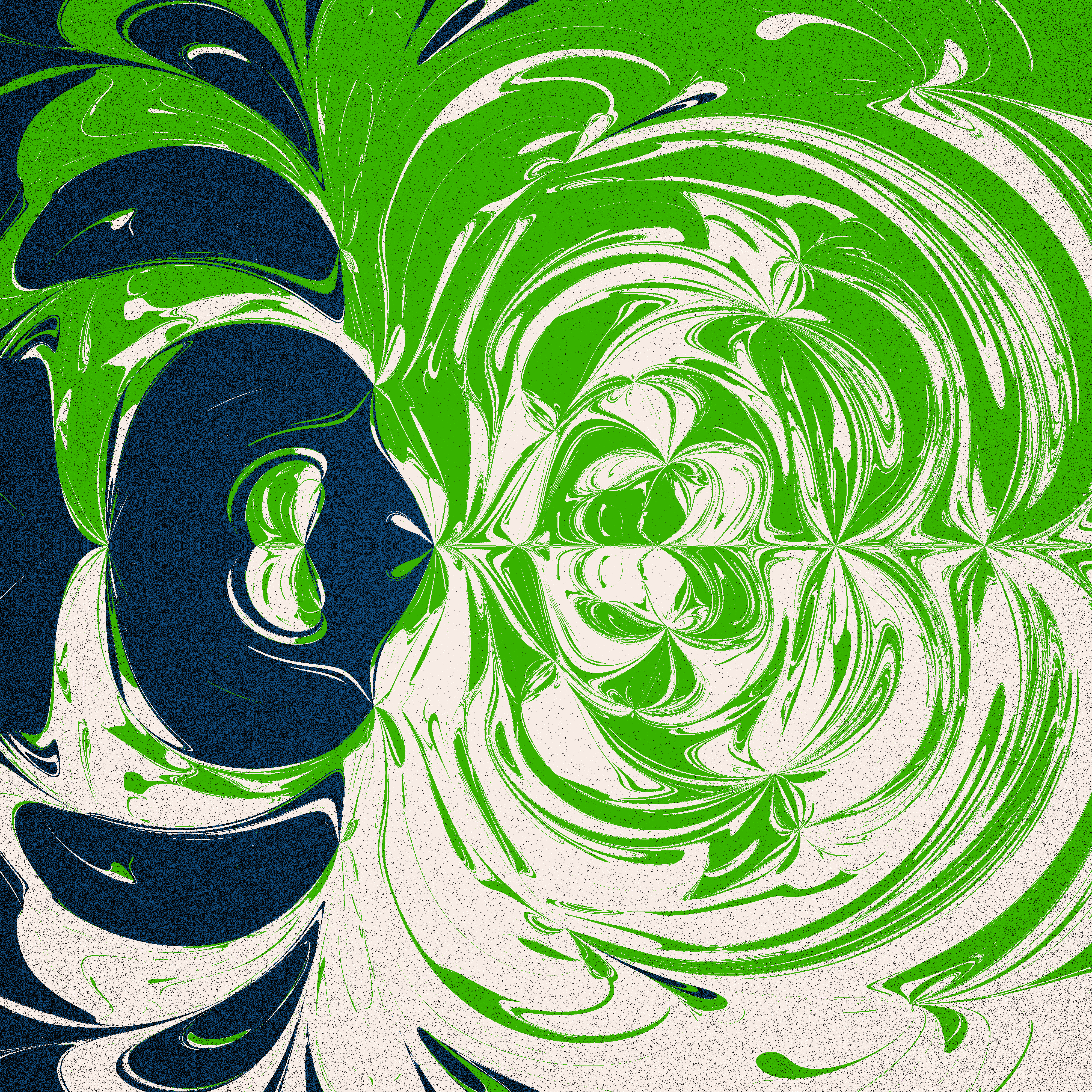}
\caption{Polynomiography of $z^3-2z+2$ under Newton (left) and Newton-Ellipsoid method.} \label{NE-Smale}
\end{figure}

\begin{figure}[h!]
\centering
\includegraphics[width=2.6in]{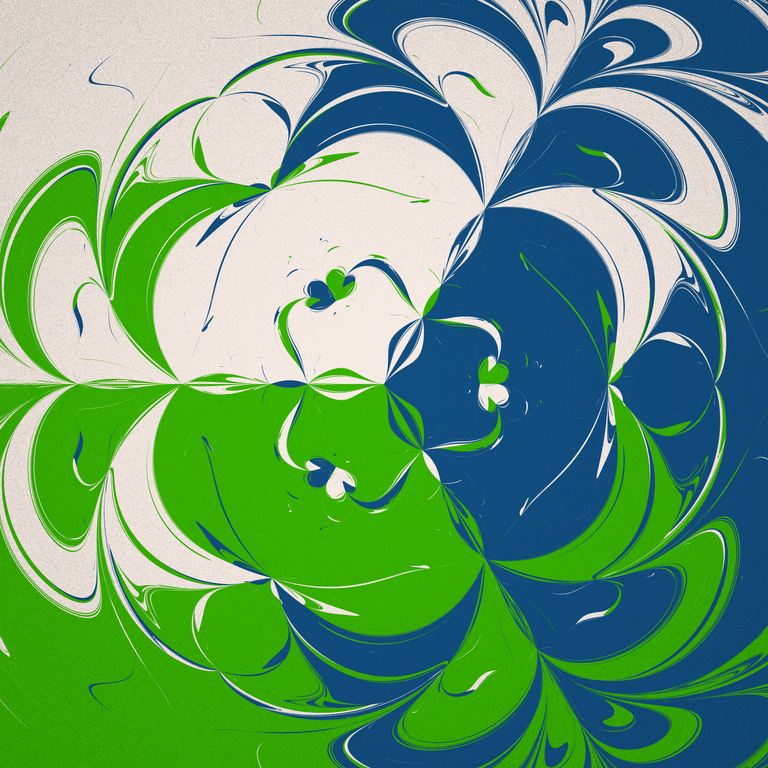}
\includegraphics[width=2.6in]{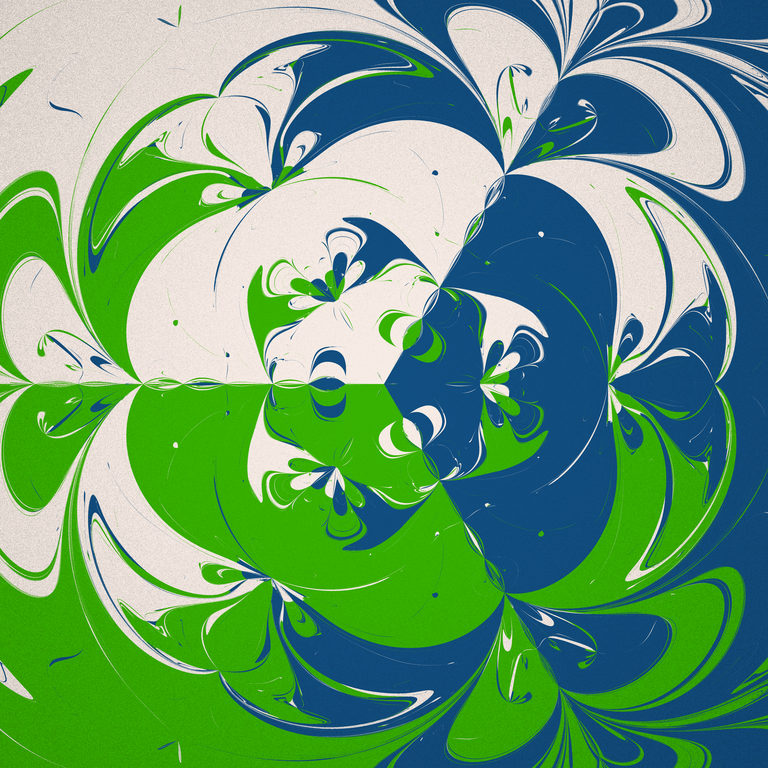}
\caption{Polynomiography of the roots of $z^3-1$ under Halley-Ellipsoid and $B_4$-Ellipsoid methods.} \label{HE-B4E}
\end{figure}

Visualization of the Newton-Ellipsoid is analogous to typical visualizations of methods such as Newton's method.  The justification in using the term polynomiography, is described in \cite{kalbook}. One typical property that is revealed via polynomiography is the basin of attractions of the roots. Starting with a grid of points around the origin, the initial rectangle $R$ containing the roots,  we took at most  60 Newton-Ellipsoid steps for each point, and then colored each point based on convergence to a root,
or lack of convergence (colored black). We examined several polynomials, $z^2-1$, $z^3-1$, $z^3-2z+2$ (the aformentioned pathological case for Newton's method). In all examples we used the rectangle with vertices $(\pm 4, \pm 4i)$. Figure \ref{FigNewton} gives the standard polynomiography of Newton's method for $z^2-1$ and $z^3-1$. Figure \ref{Newton-Ellipsoid} gives the polynomiography of Newton-Ellipsoid for the same polynomials.  Figure \ref{NE-Smale} contrasts Newton and Newton-Ellipsoid  for $z^3-2z+2$.
Figure \ref{HE-B4E} applies $B_3$-Ellipsoid and $B_4$-Ellipsoid to $z^3-1$.

\subsection{Performance of the Methods}
Looking at the image for $z^3-1$ under Newton-Ellipsoid method we notice some black dots scattered around the image. These represent points that diverge (i.e. that fail to converge to a root of $z^3-1$). We attempted to lower the number of black dots by increasing both the number of iterations and the density of the grid. However, these did not significantly seem to decrease the density of dark points. It may happen that a point $z_0$ that is near a root at the outset will move to a point far from that root in early iterations when the enclosing ellipsoid is large.

Comparing the images of $z^3-1$ and $z^3-2z+2$ under Newton's method and Newton-Ellipsoid method, we see clear differences in behavior. While using Newton's method yields more predictable convergence patterns, it also results in dense regions of divergence. However, while Newton Ellipsoid appears to be slower in speed than Newton's, it's fundamental advantage as a root-finding method is the seemingly lack of regions of divergence with positive measure. Using higher-order ellipsoid methods, with $B_3$ and $B_4$, seemed to make the general shape of basins of attraction look closer to Voronoi region of the roots. Next we describe how to improve the speed of Newton-Ellipsoid methods.

\subsection{Using Proximity Tests}

One of the important practical algorithms in computing all the roots of a polynomial is Weyl's algorithm, see \cite{pan97}.  Weyl's algorithm is two-dimensional version of the bisection algorithm. It begins with an
initial "suspect" square containing all the roots. Given a suspect square, we
partition it into four congruent subsquares. At the center of each of the
four subsquare we perform a proximity test, i.e. we estimate the distance
from the center to the nearest zero. If the proximity test guarantees that
the distance exceeds half of the length of the diagonal of the square, then the square cannot contain any zeros and it is discarded. The remaining
squares are called suspect and each of them will recursively be partitioned into four congruent subsquares and the process repeated.

One of the application of bounds stated in this article is in performing proximity tests. We refer the readers to \cite{kalbook} for details. However, we simply remark here that while using the Newton-Ellipsoid method at each input $z_0$, we can perform a proximity test and try to estimate an initial disk centered at $z_0$ that would contain a root.  This would avoid the use of large disk, hence decrease the number of iterations. We have not carried out any computational results with this approach.

\section{Concluding Remarks}
In this article we have introduced the Newton-Ellipsoid method and its generalizations.  While Newton-Ellipsoid is computationally slower than Newton's method, it compensates for this weakness by having no dense regions of divergence. Thus, it may be a practical alternative to the use of Newton's method.  Additionally, we offered some polynomiography which is interesting in its own right. The theoretical performance of the algorithm is a subject of future research.

Finally, we remark that since ellipsoid method is applicable to any dimension, assuming a generalization of Theorem \ref{FTA}, the Newton-Ellipsoid can be applied to system of polynomial equations. The utility of such method via testing is the subject of future considerations.

\bigskip


\begin{thebibliography}{9}

\bibitem{beardon91} A. F. Beardon, {\em Iteration of Rational Functions: Complex Analytic Dynamical Systems}. Springer-Verlag, New York, 1991.

\bibitem{Blum}  L. Blum, F. Cucker, M. Shub, S. Smale, \textit{Complexity and Real Computation}, Springer-Verlag, New York, 1998.\filbreak

\bibitem{cay} A. Cayley, The Newton-Fourier imaginary problem,
{\it American Journal of Mathematics}, \textbf{2} (1879) 97.\filbreak

\bibitem{kalbound} B. Kalantari, An infinite family of bounds on
zeros of analytic functions and relationship to Smale's bound, {\it
Mathematics of Computation}, \textbf{74} (2005), 841-852. \filbreak

\bibitem{kalbook} B. Kalantari, \textit{Polynomial Root-Finding and Polynomiography}, World Scientific,  Hackensack, NJ, 2008.\filbreak

\bibitem{KalDCG} B.  Kalantari,  Polynomial root-finding methods whose basins of attraction approximate Voronoi diagram,  \emph{Discrete \& Computational Geometry},  \textbf{46}  (2011) 187-203.\filbreak

\bibitem{kalFTA} B. Kalantari, A One-line proof of the Fundamental Theorem of Algebra with Newton's Method as a consequence, http://arxiv.org/abs/1409.2056,  2014. \filbreak

\bibitem{kha79} L. G. Khachiyan, A polynomial algorithm in linear programming, {\it Doklady Akademia Nauk SSSR},
(1979), 1093 - 1096.\filbreak

\bibitem{mc2005} J. M. McNamee and M. Olhovsky, A Comparison of a Priori Bounds on (Real or Complex) Roots of Polynomials, Proceedings of 17th IMACS World Congress, Scientific Computation, Applied Mathematics and
Simulation Paris, France, 2005. \filbreak

\bibitem{McMullen87} C. McMullen, Families of rational maps and iterative root-finding algorithms, {\it The Annals of Math.}, \textbf{125}  (1987) 467-493. \filbreak

\bibitem{milnor} J. Milnor, {\em Dynamics in One Complex Variable:
Introductory Lectures}, Vol 160, 3rd en. Princeton University Press, New Jersey, 2006. \filbreak

\bibitem{pan97} V. Y. Pan, Solving a polynomial equation: some
history and recent progress, \emph{SIAM Review}, \textbf{39} (1997) 187-220.\filbreak

\bibitem{peit} H.-O. Peitgen, D. Saupe, H. Jurgens, L. Yunker,
{\it Chaos and Fractals}, New York, Springer-Verlag, 1992.\filbreak

\end{thebibliography}
\end{document}